\documentclass[11pt, a4paper, reqno]{amsart}
\usepackage{preambolo}

\begin{document}


\setcounter{secnumdepth}{3}

\setcounter{tocdepth}{2}

\title[arXiv]{\textbf{Real analytic curves of almost complex structures}
}

\author[Lorenzo Sillari]{Lorenzo Sillari}

\address{Lorenzo Sillari: Scuola Internazionale Superiore di Studi Avanzati (SISSA), Via Bonomea 265, 34136 Trieste, Italy.} \email{lsillari@sissa.it}

\maketitle

\begin{abstract} 
\noindent \textsc{Abstract}. We prove that, on a compact almost complex manifold, the space of almost complex structures whose Nijenhuis tensor has rank at least $k$ at every point is either empty or dense in each path-connected component of the space of almost complex structures. In particular, this applies to maximally non-integrable almost complex structures.
\end{abstract}

\blfootnote{  \hspace{-0.55cm} 
{\scriptsize 2020 \textit{Mathematics Subject Classification}. Primary: 32Q60; Secondary: 53C15, 32G07. 
\vspace{.1cm}

\noindent The author is partially supported by GNSAGA of INdAM.}}

\section{Introduction}\label{sec:intro}

Let $M$ be a compact almost complex $2m$-manifold and let 
\[
\J = \{ J \in \End (TM) : \, J^2 = - \id \}
\]
be the space of almost complex structures on $M$. For any $J \in \J$, its Nijenhuis tensor $N_J$ is a classical invariant that encodes information on $J$. For instance, by Newlander-Nirenberg's theorem, $J$ is integrable if and only if $N_J = 0$, allowing to characterize complex structures in terms of $N_J$. Alongside complex structures, a relevant class of almost complex structures is that of \emph{maximally non-integrable} almost complex structures, i.e., those structures whose Nijenhuis tensor has maximal rank at every point of $M$. Maximally non-integrable structures represent the most generic case of almost complex structures, while complex structures the less generic one. Hence, maximally non-integrable structures are often source of many examples \cite{CGGH21, CGG23, CGHR23, CW21, CPS22}.\\
More in general, for each $k \in \N$ one can consider the space of almost complex structures whose Nijenhuis tensor has rank at least $k$ at every point, namely
\[
\J_k = \{ J \in \J: \, \, \rk N_J |_x \ge k \text{ for all } x \in M \}.
\]
Then the space of maximally non-integrable structures is given by $\J_1$ if $2m = 4$, or by $\J_m$ if $2m \ge 6$.  Almost complex structures in $\J_1$ are also called \emph{everywhere non-integrable}. Since in general $\J$ and $\J_k$ are not path-connected, we fix a path-connected component $\mathcal{C}$ of $\J$ and we consider the subsets $\mathcal{C}_k := \mathcal{C} \cap \J_k \subseteq \mathcal{C}$.\\
The problem of existence of maximally non-integrable structures has been addressed (and partially solved) by Coelho, Placini and Stelzig \cite{CPS22}. They establish an $h$-principle for the spaces $\J_k$, which has as the main consequence the fact that every almost complex structure is homotopic to a maximally non-integrable one if $2m \ge 10$ and to an everywhere non-integrable one if $2m \ge 6$ \cite[Corollary A.1]{CPS22}. In particular, maximally non-integrable (resp.\ everywhere non-integrable) structures always exist on each path-connected component $\mathcal{C}$ of $\J$ if $2m \ge 10$ (resp.\ $2m \ge 6$). Moreover, in dimension $4$ and $6$, there are necessary and sufficient topological conditions for the existence of maximally non-integrable structures (\cite[Corollaries A.2 and A.3]{CPS22}, see also \cite{Arm97, Bry06}).
\vspace{.2cm}

The main theorem of this paper is a density result for the spaces $\mathcal{C}_k$.

\begin{mythm}{\ref{thm:main}}
Let $M$ be a compact almost complex manifold and let $\mathcal{C}$ be a path-connected component of $\J$. Let $\mathcal{C}_k$ be the subspace of $\mathcal{C}$ consisting of almost complex structures of rank at least $k$ at every point of $M$. Then $\mathcal{C}_k$ is either empty or dense in $\mathcal{C}$.
\end{mythm}

As a consequence, we make precise in which sense maximally non-integrable structures are generic.

\begin{mycor}{\ref{cor:main}}
    Let $M$ be a compact almost complex manifold. Then the space of maximally non-integrable almost complex structures on $M$ is either empty or dense on each path-connected component of $\J$.
\end{mycor}

The main tools used in proving theorem \ref{thm:main} are a local result which allows to control the rank of $N_J$ in a neighborhood of a given almost complex structure (Lemma \ref{lemma:bound}), and the fact that continuous curves of almost complex structures can be approximated arbitrarily well in the $C^0$-topology of $\J$ by \emph{real analytic curves}, i.e., curves which locally behave like small deformations of almost complex structures (Lemma \ref{lemma:existence}).\\
As a consequence of the density result, almost complex structures on high-dimensional manifolds can be approximated by everywhere/maximally non-integrable ones.

\begin{mythm}{\ref{thm:approximation}}
    Let $M$ be a compact almost complex $2m$-manifold. Then every almost complex structure on $M$ can be approximated arbitrarily well in the $C^\infty$-topology by:
    \begin{itemize}
        \item an everywhere non-integrable structure if $2m \ge 6$;
        \item a maximally non-integrable structure if $2m \ge 10$.
    \end{itemize}
\end{mythm}

These results hold regardless of the dimension of $M$ for invariant structures on homogeneous spaces.
\vspace{.2cm}

In the last part of the paper, we apply the approximation result to invariants of almsot complex structures. Denote by $A^{1,0}$ the space of $(1,0)$-forms on an almost complex manifold. The number
\[
h^1_{d+d^c} := 2 \dim_\C ( \ker d \cap A^{1,0} )
\]
has been introduced in \cite{ST23b} as the dimension of the kernel of an elliptic Laplacian. In the same work, it is shown that it satisfies a topological upper bound $h^1_{d+d^c} \le b_1$ and it can distinguish between almost complex structures compatible with the same symplectic form. Applying theorem \ref{thm:approximation}, we obtain that $h^1_{d+d^c} =0$ for a generic almost complex structure on manifolds of dimension $2m \ge 10$ (Theorem \ref{thm:generically}). This result should be compared with a conjecture of Li and Zhang on the generic vanishing of the almost complex invariant $h^{-}_J$  \cite[Conjecture 2.4]{DLZ13}.\vspace{.2cm}

Overall, our results contribute to a better understanding of the space of almost complex structures of a compact almost complex manifold.

\subsection*{Acknowledgments.} The author wishes to thank the organizers of the conference ``Cohomology of Complex Manifolds and Special Structures - III" (Trento, January 2023), and the Department of Mathematics of LMU-Munich, where several ideas leading to the present paper have been discussed. He is especially thankful to Jonas Stelzig for useful conversations and comments.

\section{Curves of almost complex structures}\label{sec:curve}

In this section we study curves of almost complex structures. After briefly recalling standard material \cite{BM10, Huy05, MK06}, we discuss the notion of real analytic curve of almost complex structures, and we prove that every continuous curve can be approximated by a real analytic one
\vspace{.2cm}.

Let $M$ be a compact almost complex $2m$-manifold. For any $x \in M$, let $T_xM$ be the tangent space of $M$ at $x$. A \emph{complex structure} on $T_xM$ is and endomorphism of the tangent space $J_x \in \End (T_xM)$ such that $J^2 = -\id_{T_xM}$. The \emph{twistor bundle of $M$} (see, e.g., \cite{AHS78, OBR85}) is the fiber bundle $\tw \rightarrow M$ whose fiber $\tw_x$ is the space of complex structures on $T_x M$, i.e.,
\[
\tw_x = \{ J_x \in \End (T_xM) : J^2_x = -\id_{T_xM} \}.
\]
A smooth section of the twistor bundle corresponds to a local choice of a complex structure which depends smoothly on $x$, hence to a global almost complex structure $J$ on $M$. While the twistor bundle is defined for arbitrary even dimensional smooth manifolds, it will admit a global section if and only if $M$ is almost complex.\\
Denote by $\J$ the space of almost complex structures on $M$, or, equivalently, the space of smooth sections of $\tw \rightarrow M$. Since $\J$ is the space of sections of a smooth vector bundle, it naturally admits the structure of a Fr\'echet manifold induced by a family of seminorms. In particular, we will consider $\J$ endowed with the $C^0$-topology induced by the seminorms. Explicitly, after fixing a Riemannian metric on $M$, two points $J_0, J_1 \in \J$ are close in the $C^0$-topology if
\[
\max_{x \in M} \lVert (J_0 - J_1)|_x \rVert < \epsilon,
\]
where $\lVert \cdot \rVert$ denotes the metric induced on $\End (T_x M)$ by the choice of Riemannian metric.\\
For any point $J_0 \in \J$, the tangent space to $\J$ at $J_0$ is
\[
T_{J_0} \J = \{ L \in \End (TM): \, LJ_0 + J_0 L =0 \}.
\]
Set $I = [0, 1] \subset \R$ and denote by $\Delta_\epsilon \subset \C$ the (open) disk of radius $\epsilon$ centered at the origin.

\begin{definition}\label{def:curve}
A \emph{(real) curve of almost complex structures} is a continuous map from $I$ to $\J$. We write it as a family of almost complex structures $J_t$ depending continuously on the parameter $t \in I$. If $J_t$ depends smoothly on $t$, we say that $J_t$ is a \emph{smooth curve of almost complex structures}.
\end{definition}

\begin{definition}\label{def:deformation}
Fix $J_0 \in \J$. A \emph{small deformation} of $J_0$ is a family of almost complex structures parametrized by a complex parameter $s \in \Delta_\epsilon$ of the form
\begin{equation}\label{eq:deformation}
J_s = (I + L_s) J_0 (I + L_s)^{-1}, \quad s \in \Delta_\epsilon,
\end{equation}
where $L_s \in \End (TM)$, $L_s J_0 + J_0 L_s =0$, and $L_s = s L + o (s)$.\\
If $L_s = sL$, we say that $J_s$ is a \emph{first order deformation} of $J_0$.
\end{definition}
A small deformation of $J_0$ naturally defines a curve of almost complex structures passing through $J_0$ parametrized by $s \in \C$, $\abs{s} < \epsilon$. Taking the restriction of $s$ to the real axis, we obtain a real curve of almost complex structures parametrized by $t \in (-\epsilon, \epsilon)$, and, up to reparametrization, a curve defined in a neighborhood of $I$. If $J_1 \in \J$ is an almost complex structure close enough to $J_0$ in the $C^0$-topology of $\J$, then $J_1$ can be written as a small deformation of $J_0$.
\vspace{.2cm}

A different way of thinking of an almost complex structure $J$ is in terms of a splitting of the complexified tangent bundle
\[
TM^\C = TM^{1,0} \oplus TM^{0,1},
\]
where $TM^{1,0}$, $TM^{0,1}$ are respectively the $(\pm i)$-eigenspaces of $J$ extended to $TM^\C$. If $J_t$ is a curve of almost complex structures, it defines a family of splittings
\begin{equation}\label{eq:splitting}
TM^\C = TM^{1,0}_t \oplus TM^{0,1}_t.
\end{equation}
If $t$ is small enough and $J_t$ is induced by a small deformation of $J_0$, we can view $TM^{1,0}_t$ as a graph over $TM^{1,0}_0$, i.e.,
\[
TM^{1,0}_t = \big( \id + \Psi (t) \big) \, TM^{1,0}_0,
\]
where $\Psi(t) = \sum\limits_{j=0}^\infty \psi_j \, t^j$ is a power series with coefficients
\[
\psi_j \in T^*M^{0,1}_0 \otimes TM^{1,0}_0
\]
and $\psi_0=0$. With respect to the description of $J_t$ as a small deformation of $J_0$, one has that
\[
\Psi(t) = \frac{1}{2} ( L_t - i J_0 L_t),
\]
while $L_t$ can be recovered from $J_0$ and $J_t$ via the formula (see, e.g., \cite[Remark 3.1]{BM10})
\begin{equation}\label{eq:inverse}
    L_t = (I - J_0 J_t)^{-1} (I + J_0 J_t).
\end{equation}
Conversely, similarly to what happens in the theory of deformations of complex structures, a given power series $\Psi(t)$, with $\Psi(0)=0$, does not automatically give rise to a small deformation of almost complex structures and one has to deal with convergence issues of the power series. When we talk about small deformations of almost complex structure we always mean that the power series in $t$ is convergent in a neighborhood of $0$.\\
In order to prove the main result of this paper, we need curves of almost complex structures that are locally described by a formal deformation of almost complex structures which is actually convergent.

\begin{definition}\label{def:analytic curve}
We say that a curve of almost complex structures $J_t$ is \emph{real analytic} if for each $t_0 \in I$ there exists $\delta >0$ such that
\[
J_t = (I + L_{t-t_0}) J_{t_0} (I + L_{t-t_0})^{-1}, \quad t \in (t_0 - \delta, t_0 + \delta)
\]
and the power series described by $L_{t-t_0}$ is convergent in a neighborhood of $t_0$ with positive radius of convergence.
\end{definition}

To show existence of real analytic curves of almost complex structures, the smooth structure of $M$ is not enough and we need to use its structure of analytic manifold. We begin with some definitions.

\begin{definition}
    An \emph{analytic structure} (or \emph{$C^\omega$-structure}) on $M$ is the datum of a maximal atlas whose transition functions are analytic. $M$ endowed with a $C^\omega$-structure is an \emph{analytic manifold}. If $M$ and $N$ are analytic manifolds, an \emph{analytic function} (or \emph{$C^\omega$-function}) $f \colon M \rightarrow N$ is a function which is analytic after composing with local charts on $M$ and $N$.\\
    An \emph{analytic fiber bundle} (or \emph{$C^\omega$-fiber bundle}) over $M$ is a fiber bundle $\pi \colon E \rightarrow M$ where $E$ is an analytic manifold and $\pi$ is a surjective analytic map. An \emph{analytic section} (or \emph{$C^\omega$-section}) of a $C^\omega$-fiber bundle is a section which is also a $C^\omega$-function.\\
    A \emph{real analytic almost complex structure on $M$} is a $C^\omega$-section of the twistor bundle of $M$.
\end{definition}

A deep theorem due to Grauert and Morrey \cite{Gra58, Mor58} establishes that every smooth manifold admits a unique $C^\omega$-structure compatible with its smooth structure. Once a smooth structure is fixed, both the total spaces of its tangent bundle and twistor bundle will admit a (unique) compatible $C^\omega$-structure. Furthermore, the fiber bundle structure is compatible with the $C^\omega$-structures, endowing $\tw \rightarrow M$ with the structure of a $C^\omega$-fiber bundle.
\vspace{.2cm}

By Steenrod's approximation theorem \cite[Section 6.7]{Ste51}, every continuous section $\sigma_0$ of a smooth fiber bundle can be approximated by a smooth section $\sigma_\infty$, once a metric is fixed on the total space $B$. Furthermore, $\sigma_\infty$ can be taken in such a way that it coincides with $\sigma_0$ on a closed subset of $B$ on which $\sigma_0$ is already smooth.\\
The problem of approximating smooth sections with $C^\omega$-sections on $C^\omega$-fiber bundles was left open by Steenrod. A solution has been found by Shiga, who proved that any smooth section of a $C^\omega$-fiber bundle can be approximated by a $C^\omega$-section (see \cite[Proposition 2]{Shi64} and \cite{Shi65}). The approximation is arbitrarily good in the $C^0$-topology of the space of sections. Using Shiga's result and Grauert-Morrey's theorem, we can prove the existence of real analytic curves of almost complex structures.

\begin{lemma}\label{lemma:existence}
    Let $M$ be a compact almost complex manifold and let $J_0$, $J_1$ be two almost complex structures on $M$. Suppose that there exists a continuous curve of almost complex structures $J_t$, $t \in I$, between $J_0$ and $J_1$. Then $J_t$ can be approximated by a real analytic curve of almost complex structures.
\end{lemma}
\begin{proof}
    For an arbitrary smooth manifold $N$, continuous curves from $I$ to $N$ can be seen as continuous sections of the product fiber bundle $I \times N \rightarrow I$. Hence, we can think of $J_t$ as a section of the fiber bundle $I \times \tw \rightarrow I \times M$ which depends continuously on $t \in I$ and smoothly on $x \in M$. We want to apply Shiga's approximation theorem to the fiber bundle $I \times \tw \rightarrow I \times M$. First, endow $\tw$ and $M$ with a $C^\omega$-structure compatible with their smooth structure. Such structure always exists and it is unique by Grauert-Morrey's theorem. Hence, $J_t$ becomes a continuous section of a $C^\omega$-fiber bundle. By \cite[Proposition 2]{Shi64} we approximate it by a $C^\omega$-section $\Tilde{J}_t$. By \cite{Shi65}, the approximation can be performed keeping the endpoint fixed since $\{ 0 \} \times M$ and $\{ 1\} \times M$ are closed $C^\omega$-submanifolds of $I \times M$, hence $\Tilde{J}_0 = J_0$ and $\Tilde{J}_1 = J_1$. To prove that $\Tilde{J}_t$ is a real analytic curve of almost complex structures in the sense of definition \ref{def:analytic curve}, fix $t_0 \in I$ and let $t \in I$ be close enough to $t_0$. Consider the operator
    \[
    \Tilde{L}_t = (I - \Tilde{J}_{t_0} \Tilde{J}_t)^{-1} (I + \Tilde{J}_{t_0} \Tilde{J}_t)
    \]
    (see \eqref{eq:inverse}). Then $\Tilde{L}_t$ defines a small deformation $\Tilde{J}_t^{def}$ of $\Tilde{J}_{t_0}$. We need to prove that $\Tilde{J}_t^{def}$ coincides with the curve $\Tilde{J}_t$ in a neighborhood of ${t_0}$. Since $\Tilde{J}_t$ is analytic in $t$, so is $\Tilde{L}_t$. The same is true for $\Tilde{J}_t^{def}$, which is completely determined by $\Tilde{L}_t$ in a small neighborhood of $t_0$. Hence $\Tilde{J}_t$ and $\Tilde{J}_t^{def}$ are both analytic in $t$. Taking derivatives of \eqref{eq:deformation} and \eqref{eq:inverse}, it is immediate to check that their derivatives at $t=t_0$ coincide, so that $\Tilde{J}_t = \Tilde{J}_t^{def}$.
\end{proof}

\begin{remark}
We actually proved that a curve of almost complex structures which is analytic in $t$ is a real analytic curve in the sense of definition \ref{def:analytic curve}.
\end{remark}

\section{The space of maximally non-integrable almost complex structures}\label{sec:density}

In this section we use real analytic curves to control the rank of the Nijenhuis tensor and we prove our main results.
\vspace{.2cm}

Let $M$ be a compact almost complex manifold and let $J$ be an almost complex structure on $M$. Existence of $J$ induces a splitting
\[
A^k = \bigoplus\limits_{p+q=k} A^{p,q}
\]
of complex $k$-forms into $(p,q)$-forms. The Nijenhuis tensor of $J$ is defined on two vector fields $X$, $Y$ as
\[
N_J (X,Y) := [JX,JY] - J[JX,Y] - J[X,JY] - [X,Y].
\]

The complex rank of the complexified Nijenhuis tensor
\[
(N_J)^\C \colon TM^\C \otimes TM^\C \rightarrow TM^\C
\]
defines a map 
\begin{align*}
    \rk N_J \colon &M \longrightarrow \N \\
    &x \longmapsto \rk (N_J)^\C |_x.
\end{align*}
Equivalently, the complex rank of $N_J$ coincides with the rank of the operator
\[
\bar \mu := \pi^{0,2} \circ d \colon A^{1,0} \longrightarrow A^{0,2}.
\]
It is easy to see that for every $x \in M$, we have that $\rk N_J|_x \in \{ 0, 1 \}$ if $2m = 4$, or $\rk N_J|_x \in \{ 0, \dots, m\}$ if $2m \ge 6$. An almost complex structure is said to be \emph{maximally non-integrable} if the rank of its Nijenhuis tensor is maximum at every point $x \in M$. For brevity, we will refer to the rank of $N_J$ as the rank of $J$.\\
Let $J_t$, $t \in I$, be a smooth curve of almost complex structures. The Nijenhuis tensor $N_{J_t}$ varies smoothly with respect to $t$, and its complex rank can be computed as the rank of the operator $\bar \mu_t := \pi^{0,2}_t \circ d$, sending $A^{1,0}_t$ to $A^{0,2}_t$.
\vspace{.2cm}

Our first result allows to control the rank of $N_{J_t}$ along small deformations of almost complex structures.

\begin{lemma}\label{lemma:bound}
    Let $M$ be an almost complex $2m$-manifold. Let $J_0$ be an almost complex structure on $M$ and let 
    \[
    J_t = (I + L_t) J_0 (I + L_t)^{-1}, \quad t \in (-\epsilon, \epsilon),
    \]
    be a small deformation of $J_0$. Then for every $ a \in (0, \epsilon)$ and for every $x \in M$, we have that
    \[
    \rk N_{J_t}|_x \ge \max \{ \rk N_{J_0}|_x, \rk N_{J_{a}}|_x \}
    \]
    for all $t \in [-a, a]$ except a finite number.
\end{lemma}
\begin{proof}
    Fix $x \in M$. Let $U$ be a small neighborhood of $x$ and let $\{ \omega^j_t\}_{j=1}^m$ be a local co-frame of $(1,0)$-forms (with respect to $J_t$). The operator $\bar \mu_t$ computed on the co-frame can be written as
    \[
    \bar \mu_t \, \omega^j_t = \sum_{ k < l} G^j_{kl} (t) \,  \omega^{\bar k \bar l }_t, 
    \]
    where $G^j_{kl} \in C^\infty ( M \times (-\epsilon, \epsilon) )$ is a power series in the $t$ variable that converges on $(-\epsilon, \epsilon)$. Fix $a \in (0, \epsilon)$. Set
    \[
    k_x := \max \{ \rk N_{J_0}|_x, \rk N_{J_{a}}|_x \}.
    \]
    Let $G (t)$ be the matrix $( G^j_{kl} (t) )$ and consider all the determinants of the $k_x \times k_x$ minors of $G$. They are power series in $t$ that converge on $[-a, a]$, therefore each of them either has a finite number of zeros or identically vanishes. Since $k_x = \max \{ \rk N_{J_0}|_x, \rk N_{J_{a}}|_x \}$, at least one determinant does not identically vanish along the curve $J_t$, showing that $\rk N_{J_t}|_x \ge k_x$ for all $t \in [-a,a]$ except a finite number.
\end{proof}

\begin{remark}
A similar result valid for first order deformations on a certain class of manifolds was already proved in \cite[Theorem 4.1]{ST23a}.
\end{remark}

Let $\J$ be the space of almost complex structures on $M$. Let $\mathcal{C}$ be a path-connected component of $\J$. We set (see section \ref{sec:intro})
\[
\mathcal{C}_k := \{ J \in \mathcal{C}: \, \, \rk N_J |_x \ge k \text{ for all } x \in M \}.
\]
The results obtained in lemma \ref{lemma:bound}, together with a compactness argument, allow to control the rank of $N_J$ along real analytic curves of almost complex structures on compact almost complex manifolds.

\begin{proposition}\label{prop:global}
    Let $M$ be a compact almost complex $2m$-manifold. Let $J_0$, $J_1$ be two almost complex structures on $M$. Suppose that there exists a real analytic curve of almost complex structures $J_t$ joining $J_0$ and $J_1$. If $J_t \in \mathcal{C}_k$ for some $t \in I$, then $J_t \in \mathcal{C}_k$ for all $t \in I$ except a finite number.
\end{proposition}

\begin{proof}
    Without loss of generality, we can assume that $J_t$ is defined on $(-\epsilon, 1 + \epsilon)$. Since $J_t$ is a real analytic curve, for each $t^* \in I$ there exists a compact neighborhood $I^*$ of $t^*$ such that $J_t$ is a small deformation of $J_{t^*}$ for all $t \in I^*$. Cover $I$ with such neighborhoods and, by compactness of $I$, extract a finite subcover $\{ t_j, I_j \}_{j = 0}^r$ with $0= t_0 < t_1 < \dots < t_r = 1$ and $I_j \cap I_{j+1} \neq \emptyset$, for $j =0, \dots, r-1$. On each $I_j$, the power series describing the small deformation of $J_{t_j}$ is convergent.\\
    Since $M$ is compact, we can cover it with a finite number of neighborhoods $\{ x_s, U_s \}_{s =1}^q$ such that $T_{x_s} U_s \cong U_s \times \R^{2m}$. By assumption, there exist $\Tilde{j} \in \{0, \dots ,r\}$ and $\Tilde{t} \in I_{\Tilde{j}}$ such that $J_{\Tilde{t}} \in \mathcal{C}_k$. By lemma \ref{lemma:bound} applied on $U_s$ to the curve $J_t$, $t \in I_{\Tilde{j}}$, we have that $\rk N_{J_t}|_x \ge k$ for all $t \in I_{\Tilde{j}}$ except a finite number. Note that the result of lemma \ref{lemma:bound} applies to every point in $U_s$ for the same values of $t$. Two consecutive intervals have non-empty overlap, and we can repeat the argument on each $I_j$. We conclude that $\rk N_{J_t}|_x \ge k$ for all $x \in U_s$ and for all $t \in I$ except at most for a finite number. Finally, since we covered $M$ with a finite number of open $U_s$, we have that $J_t \in \mathcal{C}_k$ for all $t \in I$ except a finite number.
\end{proof}

As an immediate consequence, we have the following:

\begin{cor}
    Let $M$ be a compact almost complex manifold and let $J_t$, $t \in I$, be a real analytic curve of almost complex structures on $M$. Then for every $x \in M$
    \[
    \rk N_{J_t}|_x = \max_{t \in I} \{ \rk N_{J_t}|_x \}
    \]
    for all $t \in I$ except a finite number.
\end{cor}

We are ready to state and prove our main result on the density of $\mathcal{C}_k$.

\begin{theorem}\label{thm:main}
    Let $M$ be a compact almost complex manifold and let $\mathcal{C}$ be a path-connected component of $\J$. Let $\mathcal{C}_k$ be the subspace of $\mathcal{C}$ consisting of almost complex structures of rank at least $k$ at every point of $M$. Then $\mathcal{C}_k$ is either empty or dense in $\mathcal{C}$.
\end{theorem}
\begin{proof}
    Let $\mathcal{C}$ be a path-connected component of $\J$. Suppose that $\mathcal{C}_k \subseteq \mathcal{C}$ is not empty. Let $J_1 \in \mathcal{C}_k$ and let $J_0 \in \mathcal{C}$ be any almost complex structure. Let $J_t$ be a continuous curve of almost complex structures joining $J_0$ and $J_1$. Up to small perturbations, we can assume that $J_t$ is a real analytic curve by lemma \ref{lemma:existence}. By proposition \ref{prop:global} and since $J_1 \in \mathcal{C}_k$, we have that $J_t \in \mathcal{C}_k$ for all $t \in I$ except a finite number. In particular, every neighborhood of $J_0$ contains almost complex structures in $\mathcal{C}_k$.
\end{proof}

As a direct consequence, we deduce that maximally non-integrable structures are generic almost complex structures.

\begin{cor}\label{cor:main}
    Let $M$ be a compact almost complex manifold. Then the space of maximally non-integrable almost complex structures on $M$ is either empty or dense on each path-connected component of $\J$.
\end{cor}

Combining theorem \ref{thm:main} with the results of Coelho, Placini and Stelzig, we get the following approximation theorem.

\begin{theorem}\label{thm:approximation}
    Let $M$ be a compact almost complex $2m$-manifold. Then every almost complex structure on $M$ can be approximated arbitrarily well in the $C^\infty$-topology by:
    \begin{itemize}
        \item an everywhere non-integrable structure if $2m \ge 6$;
        \item a maximally non-integrable structure if $2m \ge 10$.
    \end{itemize}
\end{theorem}

\begin{proof}
    By \cite[Corollary A.1]{CPS22}, existence of everywhere (resp.\ maximally) non-integrable structures is guaranteed on each path-connected component in dimension $2m \ge 6$ (resp.\ $2m \ge 10$). Approximation in the $C^0$-topology is a direct consequence of the existence of these structures and of theorem \ref{thm:main}. To prove approximation in the $C^\infty$-topology it is enough to observe that the perturbed curve of almost complex structures $\Tilde{J}_t$ obtained by lemma \ref{lemma:existence} is actually smooth in $t$.
\end{proof}

\begin{remark}
Taking into account \cite[Corollary A.4]{CPS22}, we obtain that every invariant almost complex structure on a homogeneous manifolds can be approximated by maximally non-integrable (non-invariant) almost complex structures. Such structures represent a large source of examples where explicit computations can be performed.
\end{remark}

We conclude this section with an apllication of the approximation result which enables us to control invariants of almost complex structures in $\mathcal{J}$.\\
Let $(M,J)$ be a compact almost complex manifold. Set $d^c := J^{-1} d J$ and consider the space
\[
\mathcal{H}^1_{d+d^c} := A^1 \cap \ker d \cap \ker d^c.
\]
The space $\mathcal{H}^1_{d+d^c}$ is the kernel of a suitable $4^\text{th}$-order elliptic self-adjoint operator \cite{ST23b, ST23d, ST23c}. By compactness of $M$, its dimension
\[
h^1_{d+d^c} := \dim_\C \mathcal{H}^1_{d+d^c}
\]
is a finite number which provides a (metric-independent) invariant of $J$.

\begin{theorem}\label{thm:generically}
    Let $M$ be a compact almost complex manifold of dimension $2m \ge 10$. Then for a generic almost complex structure we have $h^1_{d+d^c} =0$.
\end{theorem}
\begin{proof}
    By theorem \ref{thm:approximation}, the space of maximally non-integrable structures is dense in $\mathcal{J}$. By \cite[Proposition 4.10]{ST23b}, if $J$ is maximally non-integrable and $2m \ge 6$ we have that $h^1_{d+d^c} (J) =0$, proving the theorem.
\end{proof}

\begin{remark}
    The results of theorem \ref{thm:generically} extend immediately to every dimension $2m \ge 6$ as long as one can show existence of maximally non-integrable structures.
\end{remark}

\printbibliography

\end{document}